# An undecidability result on limits of sparse graphs

Endre Csóka

Eötvös Loránd University, Budapest, Hungary

**Abstract**

Given a set $\mathcal{B}$ of finite rooted graphs and a radius $r$ as an input, we prove that it is undecidable to determine whether there exists a sequence $(G_i)$ of finite bounded degree graphs such that the rooted $r$-radius neighbourhood of a random node of $G_i$ is isomorphic to a rooted graph in $\mathcal{B}$ with probability tending to 1. Our proof implies a similar result for the case where the sequence $(G_i)$ is replaced by a unimodular random graph.

## 1 Introduction

There are two different approaches of handling large graphs by limits of finite graphs. The first is the Lovász–Szegedy theory[8], which is about dense graphs (with $O(|V(G_n)|^2)$ edges), and the other is about sparse graphs (with $O(|V(G_n)|)$ edges), introduced by Benjamini and Schramm[2]. Although the two concepts are developed for two different asymptotic classes of graphs, it turned out that there is a lot of analogy between the two theories, and in this paper, we show a further one.

Either model tries to describe what information we can get about a graph in constant time, using random access to the vertices. In the dense case, we can sample from the distribution of the subgraph spanned by a constant number of uniform random vertex, and in the sparse case, we can sample from the distribution of the constant radius neighbourhood of a uniform random vertex. In both cases, an important question is which the possible distributions are. Or it is better to ask what the closure of these distributions is.

Both theory use some kind of limit objects of graph sequences. In the dense case, these are symmetric measurable functions $[0,1] \times [0,1] \to [0,1]$. Contrarily, in the sparse case, a limit object is a sequence of finer and finer distributions of the neighbourhoods of the vertices with larger and larger constant radius. Therefore, describing the set of limits of all convergent graph sequences is equivauent to describing, for each $r$, the closure of the set of all distributions of the $r$-radius neighbourhoods of the vertices of the different graphs. This shows that the question has even higher importance in the sparse model than in the dense model.

About dense graphs, Hatami and Norine[7] showed the undecidability of the problem of determining the validity of a linear weak inequality between the densities. Namely, let $t_H(G)$ denote the probability that $|V(H)|$ random nodes of $G$ span $H$. Then the undecidable question is, given the sequence of graphs $H_1, H_2, ..., H_n$ and real coefficients $\lambda_1, \lambda_2, ..., \lambda_n$ as input, whether $\sum \lambda_i \cdot t_{H_i}(G) \geq 0$ holds for all graphs $G$. This result shows that the closure of the set of these distributions does not form a nice set, in this particular sense.

For sparse graphs, the closure of the set of distributions is convex, so contrary to the dense case, linear inequalities are sufficient to completely describe it. In this paper, we will show the undecidability whether this closure contains a point at which all densities of some given neighbourhoods are 0, or there exists a positive lower bound of the total relative frequency of these neighbourhoods in all graphs.

We mention that Bulitko investigated similar problems, namely he showed the undecidability of the following. Given a rooted graph $G$ as an input, does there exist a finite graph such that the 1-neighbourhood of each of its nodes is isomorphic to $G$ [4]? The main difference between the result of Bulitko and of ours is that Bulitko dealt with all neighbourhoods, while we deal with almost all neighbourhoods, in a particular sense.

We also mention that Bulitko showed that finding a finite or finding an infinite graph with the above property are not equivalent[5]. This implies the analogous nonequivalence for our case (Problem $\tilde{I}$), as well.

## 1.1 Notation and the Aldous–Lyons problem

**Graph** means **finite or infinite** graph with at least one vertex and with degrees bounded by a fix constant. **Random graph** is a probability distribution on connected rooted graphs (up to isomorphism). For a graph $G$, a node $o$ of $G$ and an integer $r$, the **$r$-neighbourhood** of $o$ is the rooted subgraph $B_r(G, o)$ spanned by all nodes of $G$ at distances at most $r$ from $o$, with the root at $o$. The $r$-neighbourhood distribution of a finite graph $G$ is the probability distribution of the $r$-neighbourhood of a uniform random node of $G$, up to isomorphism. Consider a sequence $(G_n)$ of finite graphs. If for all $r$, the distribution of the $r$-neighbourhoods of $G_n$ converges in $n$, then we say that $(G_n)$ is **convergent**. In this case, there exists a random graph with the same distribution of the $r$-neighbourhoods of the root as these limit distributions of $(G_n)$. We call this random graph the **limit** of $(G_n)$ and we call a random graph which is the limit of some sequence of finite graphs a **limit random graph**. All limit random graphs satisfy a natural condition, called **unimodularity**, which means the following. Directed-edge-rooted graph means an ordered triple $(G, o, i)$ where $G$ is a graph, and $(o, i)$ is an edge of $G$. For each random graph, we can construct a distribution of directed-edge-rooted graphs as follows. We choose a rooted graph $(G, o)$ randomly with the distribution, and we choose a uniform random neighbour $i$ of $o$. If this distribution of the $(G, o, i)$-s is invariant under the change of the direction of the edge-root (namely, the distribution is symmetric to the last two arguments) then we say that the original random graph is unimodular.

We do not have any nice description of the set of all limit random graphs. However, as Aldous and Lyons pointed out[1], there is a hope that this coincides with the set of all unimodular random graphs, which is something better understood. For example, Gábor Elek showed that each unimodular random graph can be represented by a topological graphing[6] which means a set of measure-preserving equivalence relations on a topological space with a probability measure. Whether these two sets really coincide is still an open question.

We show that describing either the set of all limit random graphs or of the unimodular random graphs is algorithmically undecidable in some sense. We also show that these sets are undecidable in the same way. Namely, we present a set of questions which are undecidable on either set, but for each question, the two answers are the same. We hope that our result can provide insight to the Aldous-Lyons problem.

## 2 Results

**Problem $L$** (limit). *Given a set $\mathcal{B}$ of finite rooted graphs and a radius $r$ as input, does there exist a limit random graph so that the $r$-neighbourhood of the root belongs a.s. to $\mathcal{B}$?*
**Problem $U$** (unimodular). *Given a set $\mathcal{B}$ of finite rooted graphs and a radius $r$ as input, does there exist a unimodular random graph so that the $r$-neighbourhood of the root belongs a.s. to $\mathcal{B}$?*
**Problem $I$** (infinite). *Given a set $\mathcal{B}$ of finite rooted graphs and a radius $r$ as input, does there exist a (possibly infinite) graph so that the $r$-neighbourhood of each node belongs to $\mathcal{B}$?*

For an input $(\mathcal{B}, r)$ denote the three answers by $L(\mathcal{B}, r)$, $U(\mathcal{B}, r)$ and $I(\mathcal{B}, r)$, respectively. We know from [2] and [1] that $L(\mathcal{B}, r) \Rightarrow U(\mathcal{B}, r) \Rightarrow I(\mathcal{B}, r)$. The first implication is derived from a kind of unimodular property of the finite graphs, while the second one holds, because if in a unimodular random graph $u$, the root a.s. has a local property, then for a randomly chosen graph from $u$, a.s. all its vertices has this property.

**Theorem 1.** *Problems $L$, $U$ and $I$ are undecidable. Moreover, there exists a subset of inputs such that, for each input $(\mathcal{B}, r)$, the three answers $L(\mathcal{B}, r)$, $U(\mathcal{B}, r)$ and $I(\mathcal{B}, r)$ are the same, and the three problems are undecidable on this set of inputs.*

**Coloured graph** means a graph with directed and coloured edges. All phenomena we defined above can also be defined for coloured graphs using a fixed finite colour set. (In this case, incidence is considered and distance in graph is measured by the undirected way.) First, we prove the same for the following two problems.

**Problem $\tilde{L}$.** *Given a set $\tilde{\mathcal{B}}$ of finite rooted coloured graphs as input, does there exist a limit random coloured graph so that the 2-neighbourhood of the root belongs a.s. to $\tilde{\mathcal{B}}$?*

**Problem $\tilde{U}$.** *Given a set $\tilde{\mathcal{B}}$ of finite rooted coloured graphs as input, does there exist a unimodular random coloured graph so that the 2-neighbourhood of the root belongs a.s. to $\tilde{\mathcal{B}}$?*

**Problem $\tilde{I}$.** *Given a set $\tilde{\mathcal{B}}$ of finite rooted coloured graphs as input, does there exist a (possibly infinite) graph so that the 2-neighbourhood of each node belongs to $\tilde{\mathcal{B}}$?*

For an input $\tilde{\mathcal{B}}$, denote the three answers by $L(\tilde{\mathcal{B}})$, $U(\tilde{\mathcal{B}})$ and $I(\tilde{\mathcal{B}})$, respectively. $L(\tilde{\mathcal{B}}) \Rightarrow U(\tilde{\mathcal{B}}) \Rightarrow I(\tilde{\mathcal{B}})$ holds by the same reason as in the uncoloured case.

**Proposition 2.** *There exists a subset of inputs such that, for each input $\tilde{\mathcal{B}}$, the three answers $L(\tilde{\mathcal{B}})$, $U(\tilde{\mathcal{B}})$ and $I(\tilde{\mathcal{B}})$ are the same, and the three problems $\tilde{L}$, $\tilde{U}$ and $\tilde{I}$ are undecidable on this set of inputs.*

*Proof of Proposition 2.* We will show the undecidablity of these problems even if $\tilde{\mathcal{B}}$ contains graphs only with degrees bounded by 4.

We prove this undecidability by reduction from Wang's Tiling Problem[9] which is the following. We get a finite set of equal-sized square tiles with a color on each edge. We want to place one of these on each square of a regular square grid so that abutting edges of adjacent tiles have the same color. The tiles cannot be rotated or reflected, but each tile can be used arbitrary many times, or in other words, we have infinitely many tiles of each kind in the set. The question is, given the set of tiles as input, whether there exists a tiling of the whole plane. Berger showed in 1966 that this question is algoritmically undecidable[3].

For each specific input of the tiling problem, we construct a set $\tilde{\mathcal{B}}$ such that the existence of a tiling of the plane is equivalent to all of $L(\tilde{\mathcal{B}})$, $U(\tilde{\mathcal{B}})$ and $I(\tilde{\mathcal{B}})$. This will imply the undecidabilities.

To each tiling of a part of the plane, we assign a graph, as follows. The nodes of this graph represent the different tiles. The colour set of the edges is $\{\rightarrow, \uparrow\} \times$ (all colours in the tiles). If a tile $Y$ is the right or up neighbour of $X$, and the colour of the abutting edges is $c$ then we put a directed edge from $X$ to $Y$ with colour $(\rightarrow, c)$ or $(\uparrow, c)$, respectively.

Consider graphs corresponding to all $5 \times 5$ tilings. We define $\tilde{\mathcal{B}}$ as the set of 2-neighbourhoods of the central nodes of all such graphs. Figure 1/b illustrates such a neighbourhood, corresponding to the tiling in Figure 1/a.

Now we prove the equivalence of $L(\tilde{\mathcal{B}})$, $U(\tilde{\mathcal{B}})$, $I(\tilde{\mathcal{B}})$ and the existence of a tiling of the plane, by the following implications.

- $L(\tilde{\mathcal{B}}) \Rightarrow U(\tilde{\mathcal{B}}) \Rightarrow I(\tilde{\mathcal{B}})$, as we have already seen.

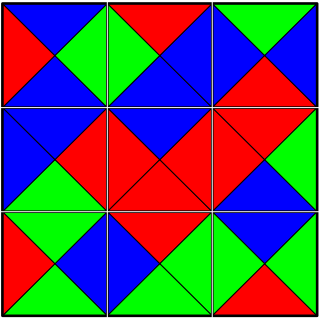 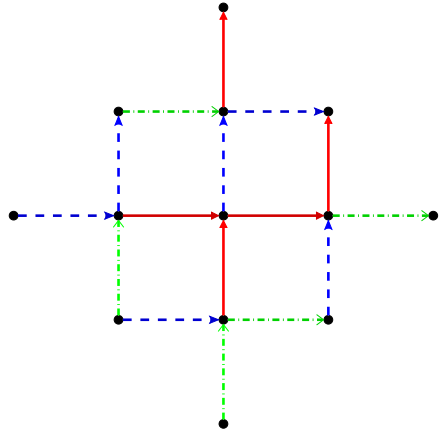 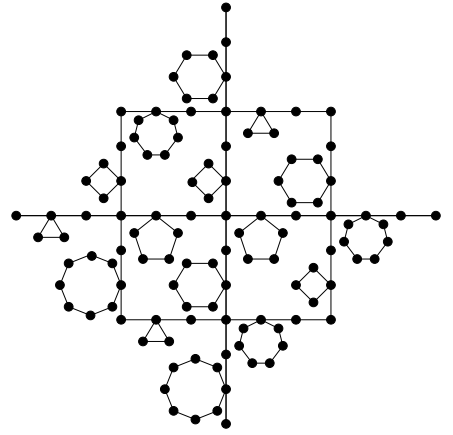

Figure 1/a     Figure 1/b     Figure 1/c

- If $I(\tilde{\mathcal{B}})$, then there exists a tiling of the plane.
  *Proof.* Consider such a graph. Each of its vertices can be represented by a tile, namely the four colours of the tile are the four colours of the incident edges in the appropriate directions. Consider the lattice group, which is the group generated by two elements, called *up* and *right*, with the defining relation $up \cdot right = right \cdot up$. This group acts on the graph, because the structure of the neighbourhoods guarantees that for each vertex $v$, $up\big(up^{-1}(v)\big) = right\big(right^{-1}(v)\big) = v$, and $up\big(right(v)\big) = right\big(up(v)\big)$. Let us take an arbitrary vertex $v$. Let us naturally identify the unit squares of the grid with the group elements. Then let the tile at each square $s$ be the representing tile of $s(v)$. This is clearly a valid tiling.

- If there exists a tiling of the plane, then $L(\tilde{\mathcal{B}})$.
  *Proof.* Consider a tiling and let $G_n$ be the graph corresponding to an $n \times n$ tiling. For each node in distance at least 2 from the border of the square, the 2-neighbourhood of it is in $\tilde{\mathcal{B}}$. Thus the relative frequency of these nodes tends to 1, so the limit of the sequence $(G_n)$ satisfies the requirements. $\square$

*Proof of Theorem 1.* We construct a recursive function $F$ that maps every finite set $\tilde{\mathcal{B}}$ of finite rooted coloured graphs to a pair $(\mathcal{B}, r)$ so that $L(\mathcal{B}, r)$, $U(\mathcal{B}, r)$, $I(\mathcal{B}, r)$, $L(\tilde{\mathcal{B}})$, $U(\tilde{\mathcal{B}})$, $I(\tilde{\mathcal{B}})$ are all equivalent. This will prove the theorem, because if there were an algorithm computing $L(\mathcal{B}, r)$, then we could compute $L(\tilde{\mathcal{B}})$, as well, by transforming $\tilde{\mathcal{B}}$ to $(\mathcal{B}, r)$ and then computing $L(\mathcal{B}, r)$.

For the sake of simplicity, we assume that the colours in $\tilde{\mathcal{B}}$ are $\{1, 2, ..., k\}$. Take a 3-length path as a graph, let $A$, $B$, $C$ and $D$ be the four consecutive nodes of it. Take an $l$-length cycle, and identify one node of the cycle with $B$. If $l = 2c+1$ or $l = 2c+2$, then let us call this graph $R_c$ and $U_c$, respectively. From a coloured graph $\tilde{G}$, we construct a graph $G = f(\tilde{G})$ with degrees bounded by 4, like from Figure 1/b to Figure 1/c, as follows. Instead of each directed edge from $X$ to $Y$ with colour $(\rightarrow, c)$ or $(\uparrow, c)$, we add a new copy of $R_c$ or $U_c$, respectively, disjointly to the graph, and identify $X$ with $A$ and $Y$ with $D$.

Let $\mathcal{G}$, $\tilde{\mathcal{G}}$ and $\tilde{\mathcal{G}}_k$ denote the isomorphism classes of graphs, coloured graphs and coloured graphs with colour set $\{1, 2, ..., k\}$, respectively. Let $\mathcal{G}_*$, $\tilde{\mathcal{G}}_*$ and $(\tilde{\mathcal{G}}_k)_*$ denote their rooted versions.

We will use a global retransformation $f^{-1} : \mathcal{G} \rightarrow \tilde{\mathcal{G}} \cup \{error\}$ and for each $k \in \mathbb{N}$, a local retransformation algorithm $f_k^{-1} : \mathcal{G}_* \rightarrow (\tilde{\mathcal{G}}_k)_* \cup \{error\}$, satisfying the following properties. (Note that $f^{-1}$ denotes *not* the inverse of $f$, but something very similar.)

We define $f^{-1}(G)$ as follows. We call a node **central** if it has degree 4 and it has two neighbours with degree 4. These nodes will be the vertices of $f^{-1}(G)$. In $G$, let us cut these nodes into 4 nodes with degrees 1, keeping all edges. Then each component should be an $R_i$ or $U_i$ – otherwise we output an *error* – and we add the corresponding directed coloured edge between the corresponding nodes of $f^{-1}(G)$. The coloured graph we get is $f^{-1}(G)$.

We define $f_k^{-1}(G_*) = f_k^{-1}(G, o)$ as follows. We find the central node $\tilde{o}$ closest to $o$. This should be unique and in distance at most $k+2$ from $o$. $\tilde{o}$ will be the root of the coloured graph. We take all central nodes in the 6-neighbourhood of $\tilde{o}$ as the vertices of the coloured graph. We cut these nodes into 4 nodes with degrees 1, keeping all edges. We should get components with diameter at most $k+3$. We decode them into coloured directed edges between the nodes. If everything was right then we output the 2-neighbourhood of the root in this coloured graph, otherwise we output an *error*.

These retransformation algorithms have the following properties.

1. $f^{-1}\big(f(\tilde{G})\big) = \tilde{G}$

2. $f_k^{-1}$ depends only on the constant radius neighbourhood of the root, namely
$$f_k^{-1}(G_*) = f_k^{-1}\big(B_{2k+8}(G_*)\big)$$

3. If $f_k^{-1}(G_*) \in (\tilde{\mathcal{G}}_k)_*$, then each vertex of $f_k^{-1}(G_*)$ is in distance at most 2 from the root.

4. A graph is globally retransformable if and only if it is locally retransformable at each point. Namely,
$$\big(\forall o \in V(G) : f_k^{-1}(G, o) \in (\tilde{\mathcal{G}}_k)_*\big) \Leftrightarrow \big(f^{-1}(G) \in \tilde{\mathcal{G}}_k\big).$$

5. There exists a constant $C_k$ that if $\tilde{G} = f^{-1}(G) \in \tilde{\mathcal{G}}_k$, then there exists a mapping $m : V(G) \to V(\tilde{G})$ satisfying the followings.
$$\forall v \in V(G) : \ B_2\big(m(v)\big) \cong f_k^{-1}(G, v)$$
$$\forall \tilde{v} \in V(\tilde{G}) : \ 1 \le \big|m^{-1}(\tilde{v})\big| \le C_k$$

Properties 1, 2, 3 and 4 are obvious consequences of the retransformating methods. At property 5, $m(v)$ can be constructed as follows. We take the central node $v_0$ closest to $v$. In the definition of $f^{-1}$, the new nodes are constructed from the central nodes of the original graph. Let $m(v)$ be the node of $\tilde{G}$ corresponding to $v_0$. Then $B_2\big(m(v)\big) \cong f_k^{-1}(G, v)$, by the definition of $f_k^{-1}$. About $m^{-1}(\tilde{v})$, it contains the central node $v_0$ of $G$ corresponding to $\tilde{v}$, and each node in $m^{-1}(\tilde{v})$ must be in distance at most $k+2$ from $v_0$, according to the construction of $m$. That is why, $\big|m^{-1}(\tilde{v})\big| \le 1 + 4 + 4^2 + ... + 4^{k+2} = C_k$.

Let $\mathcal{B} = F(\tilde{\mathcal{B}})$ be the set of rooted graphs $G_*$ satisfying that its degrees are at most 4, and its nodes are in distance at most $r = 2k + 8$ from the root, and $f_k^{-1}(G_*) \in \tilde{\mathcal{B}}$. Clearly, $F$ is a recursive function.

If $L(\tilde{\mathcal{B}})$, then there exists a finite coloured graph sequence $(\tilde{G}_n)$ so that the relative frequency of the nodes with 2-neighbourhoods from $\tilde{\mathcal{B}}$, tends to 1. Then using Property 5, we get that the relative frequency of $r$-neighbourhoods of $f(\tilde{G}_n)$ which are from $\mathcal{B}$, tends to 1. This implies $L(\mathcal{B}, r)$.

If $I(\mathcal{B}, r)$, then take a graph $G$ so that $\forall v \in V(G) : B_r(v) \in \mathcal{B}$. By Property and 4, $f^{-1}(G) \in \tilde{\mathcal{G}}_k$, and $\forall v \in V(G) : f_r^{-1}(G, v) \in \tilde{\mathcal{B}}$. Thus, $f^{-1}(G)$ shows that $I(\tilde{\mathcal{B}})$.

To sum up, $L(\mathcal{B}, r) \Rightarrow U(\mathcal{B}, r) \Rightarrow I(\mathcal{B}, r) \Rightarrow I(\tilde{\mathcal{B}}) \Leftrightarrow U(\tilde{\mathcal{B}}) \Leftrightarrow L(\tilde{\mathcal{B}}) \Rightarrow L(\mathcal{B}, r)$, which proves the equivalences. □

# 3 Acknowledgements

The author thanks to Gábor Elek for his suggestion to think on, and his observation that the proof works not only on limit random graphs but on unimodular random graphs, as well.
Thanks to László Lovász for his help making this paper.